\input amstex.tex
\input amsppt.sty
\magnification=\magstep1
\baselineskip=16truept
\vsize=22.2truecm
\NoBlackBoxes
\nologo
\pageno=1
\topmatter
\def\Z{\Bbb Z}

\def\F{\Bbb F}
\def\l{\left}
\def\r{\right}
\def\bg{\bigg}
\def\({\bg(}
\def\[{\bg[}
\def\){\bg)}
\def\]{\bg]}
\def\t{\text}
\def\f{\frac}

\def\em{\emptyset}
\def\se {\subseteq}

\def\sm{\setminus}

\def\cs{\ldots}

\def\al{\alpha}

\def\da{\delta}

\def\colon{{:}\;}
\def\Proof{\noindent{\it Proof}}
\def\Remark{\medskip\noindent{\it  Remark}}
\def\Ack{\medskip\noindent {\bf Acknowledgment}}
\hbox {J. Combin. Theory Ser. A 100(2002), no. 2, 387--393.}
\medskip
\title {A lower bound for $|\{a+b\colon a\in A,\ b\in B,\ P(a,b)\not=0\}|$}\endtitle
\author Hao Pan and Zhi-Wei Sun\endauthor
\medskip
\abstract
Let $A$ and $B$ be two finite subsets of a field $\F$.
In this paper we provide a nontrivial lower bound for
$|\{a+b\colon a\in A,\ b\in B,\ \t{and}\ P(a,b)\not=0\}|$
where $P(x,y)\in \F[x,y]$.
\endabstract
\thanks 2000 {\it Mathematics Subject Classification}.
Primary 11B75; Secondary 05A05, 11C08.
\newline\indent The second author is responsible for all the communications,
and supported by
the Teaching and Research Award Program for Outstanding Young Teachers
in Higher Education Institutions of MOE, and
the National Natural Science Foundation of P. R. China.
\endthanks
\endtopmatter
\document

\heading{1. Introduction}\endheading

Let $\F$ be a field and let $\F^{\times}$ be the multiplicative group $\F\setminus\{0\}$.
The additive order of the (multiplicative) identity of $\F$ is
either infinite or a prime, we call it the
{\it characteristic} of $\F$.

Let $A$ and $B$ be finite subsets of the field $\F$. Set
$$A+B=\{a+b\colon a\in A\ \t{and}\ b\in B\}$$
and $$ A\dotplus B=\{a+b\colon a\in A,\ b\in B,\ \t{and}\ a\not=b\}.$$
The theorem of Cauchy and Davenport (see, e.g., [N, Theorem 2.2]) asserts that
if $\F$ is the field of residues modulo a prime $p$, then
$$|A+B|\geq\min\{p,\ |A|+|B|-1\}.$$
In 1964 Erd\H os and Heilbronn (cf. [EH] and [G]) conjectured that
in this case
$$|A\dotplus A|\geq\min\{p,\ 2|A|-3\},$$
this was confirmed by Dias da Silva and
Hamidoune [DH] in 1994.
In 1995--1996  Alon, Nathanson and Ruzsa [ANR1, ANR2] proposed
a polynomial method to handle similar problems, they showed that
if $|A|>|B|>0$ then
$$|A\dotplus B|\geq\min\{p,\ |A|+|B|-2\}$$
where $p$ is the characteristic of the field $\F$.
The method usually yields a nontrivial conclusion provided that certain
coefficient of a polynomial, related in some special way to the
additive problem under considerations, does not vanish.

What can we say about the
cardinality of the restricted sumset
$$C=\{a+b\colon a\in A,\ b\in B,\ \t{and}\ P(a,b)\not=0\}\tag1$$
where $P(x,y)\in\F[x,y]$?
We will make progress in this direction
by relaxing (to some extent) the limitations of the polynomial
method. Our approach allows one to draw conclusions even if
no coefficients in question are known explicitly.

Throughout this paper, for $k,l\in\Z$ each of the intervals $(k,l),[k,l),(k,l],[k,l]$
will represent the set of {\it integers} in it.
For a polynomial
$P(x_1,\cs,x_n)$ over a field, we let $\hat P(i_1,\cs,i_n)$
stand for the coefficient of $x_1^{i_1}\cdots x_n^{i_n}$ in
$P(x_1,\cs,x_n)$.

Let $\Bbb E$
be an algebraically closed field and $P(x)$ be a polynomial over $\Bbb E$.
For $\al\in\Bbb E$, if
$(x-\al)^m\mid P(x)$ but $(x-\al)^{m+1}\nmid P(x)$, then
we call $m$ the {\it multiplicity} of $\al$ with respect to $P(x)$
and denote it by $m_P(\al)$. For any positive integer $q$, we set
$$N_q(P)=q|\{\al\in\Bbb E^{\times}\colon  m_P(\al)\ge q\}|
-\sum_{\al\in\Bbb E^{\times}}\{m_P(\al)\}_q\tag2$$
where $\{m\}_q$ denotes the least nonnegative residue of $m\in\Z$ modulo $q$.
Note that $N_1(P)$ is the number of distinct roots in $\Bbb E^{\times}$
of the equation $P(x)=0$. Let $p$ be the characteristic of $\Bbb E$, and
$$\Cal P(p)=\cases\{1,p,p^2,\cs\}&\t{if}\
p<\infty,\\\{1\}&\t{otherwise}.\endcases$$
We also define
$$N(P)=\max_{q\in \Cal P(p)}q|\{\al\in\Bbb E^{\times}\sm\{-1\}\colon m_P(\al)\ge q\}|.
\tag3$$
Clearly $N(P)\le \sum_{\al\in \Bbb E^{\times}\sm\{-1\}}m_P(\al)\le\deg P(x)$.

Let $\F$ be a field of characteristic $p$, and let $\Bbb E$
be the algebraic closure of $\F$. Any $P(x)\in\F[x]$ can be viewed
as a polynomial over $\Bbb E$ so that $N_q(P)\ (q=1,2,3,\cs)$ and
$N(P)$ are well defined.
If $P(x)\in\F[x]$ is irreducible and it
has a repeated zero in $\Bbb E$, then $p<\infty$ and $P(x)=f(x^p)$
for some irreducible $f(x)\in \F[x]$ (see, e.g. [W, Theorem 9.7]);
as $x^p-\al^p=(x-\al)^p$ for all $\al\in\Bbb E$, by induction we find
that the multiplicity of any zero of $P(x)$ belongs to $\Cal P(p)$.

The key lemma of this paper is the
following new result.

\proclaim{Lemma 1} Let $P(x)$ be a polynomial over the field $\F$
of characteristic $p$.
Suppose that there exist nonnegative
integers $k<l$ such that $\hat P(i)=0$ for all $i\in(k,l)$. Then
either $x^l\mid P(x)$, or $\deg P(x)\le k$, or $N_q(P)\ge l-k$
for some $q\in\Cal P(p)$.
\endproclaim

With helps of Lemma 1 and the polynomial method, we are able to obtain
the following main result.

 \proclaim{Theorem 1} Let $\F$ be a field of characteristic $p$,
and let $A$ and $B$ be two finite nonempty subsets of $\F$.
Furthermore, let $P(x,y)$ be a polynomial over $\F$ of degree
$d=\deg P(x,y)$ such that for some $i\in[0,|A|-1]$
and $j\in[0,|B|-1]$ we have $\hat P(i,d-i)\not=0$ and $\hat
P(d-j,j)\not=0$.
Define $P_0(x,y)$ to be the homogeneous polynomial of degree $d$
such that $P(x,y)=P_0(x,y)+R(x,y)$ for some $R(x,y)\in\F[x,y]$
with $\deg R(x,y)<d$, and put $P^*(x)=P_0(x,1)$.
Then, for the set $C$ given by $(1)$, we have
$$|C|\geq \min\{p-m_{P^*}(-1),\ |A|+|B|-1-d-N(P^*)\}.\tag4$$
\endproclaim

\Remark\ 1. In the case $d=\deg P(x,y)=0$, Theorem 1 yields the
Cauchy-Davenport theorem.

\smallskip

Lemma 1 and Theorem 1 will be proved in the next section.

Now we give some consequences of Theorem 1.

\proclaim{Corollary 1} Let $\F$ be a field of characteristic $p$,
and let $A$ and $B$ be finite subsets of $\F$. Let $k,m,n$ be nonnegative integers and
$Q(x,y)\in\F[x,y]$ have degree less than $k+m+n$. If $|A|>k$ and
$|B|>m$, then
$$\aligned
&|\{a+b\colon a\in A,\ b\in B,\ \t{and}\ a^kb^m(a+b)^n\not=Q(a,b)\}|
\\&\ \geq\min\{p-n,\ |A|+|B|-k-m-n-1\}.\\
\endaligned\tag 5$$
\endproclaim
\Proof.
For $P(x,y)=x^ky^m(x+y)^n-Q(x,y)$, clearly $\hat P(k,m+n)=\hat P(k+n,m)=1$
and $P^*(x)=x^k(x+1)^n$.
Since $N(P^*)=0$, the desired result follows from Theorem 1. \qed

\Remark\ 2. When $k=m=1,\ n=0$ and $Q(x,y)=1$, our Corollary 1
yields [ANR1, Theorem 4] which is also [ANR2, Proposition 4.1].

\proclaim{Corollary 2} Let $\F$ be a field of characteristic $p\not=2$,
and let $A,B$ and $S$ be finite nonempty
subsets of $\F$. Then
$$|\{a+b\colon a\in A,\ b\in B,\ \t{and}\ a-b\not\in
S\}|\ge\min\{p,\ |A|+|B|-|S|-q-1\}\tag6$$
where $q$ is the largest element of $\Cal P(p)$ not exceeding $|S|$.
\endproclaim
\Proof. Let $C=\{a+b\colon a\in A,\ b\in B,\ \t{and}\ a-b\not\in S\}$.
By applying Theorem 1 with $P(x,y)=\prod_{s\in S}(x-y-s)$,
we obtain the desired lower bound for $|C|$.
\qed

\Remark\ 3. In the case $S=\{0\}$, Corollary 2 was first obtained
by Alon, Nathanson and Ruzsa [ANR1, ANR2]. When $|A|=|B|=k$, $2\mid |S|$ and $|S|<p$,
the lower bound in (6) can be replaced by
$\min\{p,2k-|S|-1\}$ as pointed out by Hou and Sun
[HS]. For a field $\Bbb F$ with $|\Bbb F|=2^n>2$, if $A,S\se\Bbb F$,
$|A|>2^{n-1}+1$ and $|S|=2^n-1$, then
$|\{a+b\colon a\in A,\ b\in \F,\ \t{and}\ a-b\not\in S\}|
=|(A+\F)\setminus S|=|\F\sm S|=1<\min\{2,|A|+|\F|-|S|-2^{n-1}-1\}$. So we cannot omit
the condition $p\not=2$ from Corollary 2.

\proclaim{Corollary 3} Let $\F$ be a field of characteristic $p$,
and let $A$ and $B$ be finite nonempty subsets of
$\F$. Let $\em\not= S\se\F^{\times}\times\F$ and $|S|<\infty$.
Then
$$\aligned& |\{a+b\colon a\in A,\ b\in B,\
\t{and}\ a+ub\not=v\ \t{if}\ \langle u,v\rangle\in S\}|
\\\ge&\min\{p-|\{v\in\F\colon\langle 1,v\rangle\in S\}|,\ |A|+|B|-2|S|-1\}.
\endaligned\tag7$$
\endproclaim
\Proof.  Just apply Theorem 1 with $P(x,y)=\prod_{\langle u,v\rangle\in
S}(x+uy-v)$ and note that $N(P^*)\le\deg P^*=|S|$. \qed

\Remark\ 4. When $p=\infty$, Corollary 3 is essentially [S, Theorem 1.1]
in the case $n=2$.

\heading{2. Proofs of Lemma 1 and Theorem 1}\endheading

\noindent{\it Proof of Lemma 1}.
We use induction on $\deg P(x)$.
When $P(x)$ is a constant, we need do nothing.
So we let $\deg P(x)>0$ and proceed to the induction step.

Write $P(x)=x^hQ(x)$ where $h=m_P(0)$ and $Q(x)\in \F[x]$. If $h<l$,
then $h\le k$ since $\hat P(i)=0$ for any $i\in (k,l)$,
therefore $\hat Q(j)=0$ for all $j\in(k-h,l-h)$.
So, without loss of generality, it can be assumed that
$P(0)\not=0$ and that $P(x)$ is monic.

Let $\Bbb E$ be the algebraic closure of the field $\F$.
Write $P(x)=\prod_{j=1}^n(x-\al_j)^{m_j}$
where $\al_1,\cs,\al_n$ are distinct elements of $\Bbb E^{\times}$
and $m_1,\cs,m_n$ are positive integers. For
$j=1,\cs,n$ let $P_j(x)=P(x)/(x-\al_j)$. As
$P(x)=P_j(x)(x-\al_j)$, $\hat P(i+1)=\hat P_j(i)-\al_j\hat
P_j(i+1)$ for $i=0,1,2,\cs.$ Note that $\hat P_j(i)=\al_j\hat
P_j(i+1)$ for every $i\in[k,l-1)$. Therefore
$$\hat P_j(i)=\al_j^{l-1-i}\hat P_j(l-1)\ \ \ \t{for all}\ i\in[k,l).
\tag 8$$
Since $P'(x)=\sum_{j=1}^n m_jP_j(x)$, we have
$$\sum_{j=1}^n m_j\hat P_j(i)=0\ \ \ \t{for any}\ i\in[k,l-1).\tag9$$
Combining (8) and (9) we find that
$$\sum_{j=1}^n m_j\al_j^{l-1-i}\hat P_j(l-1)=0\ \ \ \t{for each}\ i\in
[k,l-1).\tag10$$

Suppose that $N_q(P)<l-k$ for any $q\in\Cal P(p)$. Then $n=N_1(P)\le l-1-k$,
hence by (10) we have
$$\sum_{j=1}^n\al_j^{s}(m_j\hat P_j(l-1))=0\ \ \ \t{for every}\
s=1,\cs,n.$$
Since the Vandermonde determinant $\|\al_j^{s}\|_{1\le s,j\le n}$
does not vanish, by Cramer's rule we have
$m_j\hat P_j(l-1)=0$ for all $j=1,\cs,n.$
Thus, in light of (8), $m_j\hat P_j(i)=0$ for any $i\in[k,l)$ and $j\in[1,n]$.

{\tt Case 1}. $p=\infty$, or $p\nmid m_j$ for some $j\in[1,n]$.

In this case there is a $j\in[1,n]$ such that $\hat P_j(i)=0$ for all $i\in(k-1,l)$.
Clearly $k>0$ since $\hat P_j(0)=P_j(0)\not=0$. Also
$N_1(P_j)\le n=N_1(P)$, and $N_q(P_j)=N_q(P)+1$ if $p<\infty$ and
$q\in\Cal P(p)\setminus\{1\}$. Thus
$N_q(P_j)\le N_q(P)+1\le l-k<l-(k-1)$ for all $q\in\Cal P(p)$.
In view of the induction hypothesis, we
should have $\deg P_j\le k-1$ and hence $\deg P(x)\le k$.

{\tt Case 2}. $p<\infty$, and $p\mid m_j$ for all $j\in[1,n]$.

In this case, $T(x)=\prod_{j=1}^n(x-\al_j)^{m_j/p}\in\Bbb E[x]$
and therefore
$P(x)=T(x)^p=(\sum_{i\ge0}\hat T(i)x^i)^p=\sum_{i\ge0}\hat T(i)^px^{ip}.$
For any real number $r$ let $\lfloor r\rfloor$ denote
the greatest integer not exceeding $r$.
Then $\lfloor k/p\rfloor\le\lfloor(l-1)/p\rfloor$ since $k\le l-1$.
Whenever $i\in(\lfloor
k/p\rfloor,\lfloor(l-1)/p\rfloor]$, we have $k<ip<l$
and hence $\hat T(i)^p=\hat P(ip)=0$.

If $q\in\Cal P(p)$ then
$$N_q(T)=\f{N_{pq}(P)}p\le\f{l-k-1}p<\bg(1+\l\lfloor\f{l-1}p\r\rfloor\bg)
-\l\lfloor\f kp\r\rfloor.$$
By the induction hypothesis, $\deg T\le\lfloor
k/p\rfloor$ and hence $\deg P=p\deg T\le k$.

So far we have completed the induction proof. \qed.

\medskip
\noindent{\it Proof of Theorem 1}. Set $k_1=|A|-1$ and $k_2=|B|-1$.
Clearly (4) holds if $|C|\ge k_1+k_2-d+1$. So we assume that
$|C|\le k_1+k_2-d$ and let $\da=k_1+k_2-d-|C|$.

Since $\hat P(d-j,j)\not=0$ for some $j\in[0,k_2]$,
$Q(x,y)=P(x,y)/\prod_{b\in B}(y-b)\not\in \F[x,y]$
(otherwise $\hat P(d-j,j)$ is zero because it equals the
coefficient of $x^{d-j}y^j$ in $y^{|B|}Q(x,y)$).
Thus there exists a $b_0\in B$ such that $P(x,b_0)$ does not
vanish identically, hence $P(a,b_0)=0$ for at most $d$ elements
$a\in\F$.
Therefore $$|C|\ge|\{a+b_0\colon a\in A\ \t{and}\
P(a,b_0)\not=0\}|\ge|A|-d$$
and so $\da<k_2$. Similarly we have $\da<k_1$.

Put $$f(x,y)=P(x,y)\prod_{c\in C}(x+y-c)\ \ \t{and}\ \ f_0(x,y)=P_0(x,y)(x+y)^{|C|}.$$
Clearly $\deg f(x,y)=\deg f_0(x,y)=d+|C|=k_1+k_2-\da$.
Let $\kappa_1\in[k_1-\da,k_1]$. Then
$\kappa_2=k_1+k_2-\da-\kappa_1\in(0,k_2]$.
As $\kappa_1+\kappa_2=\deg f(x,y)$ and $f(x,y)$ vanishes over the Cartesian product
$A\times B$, $\hat f(\kappa_1,\kappa_2)=0$ by [A, Theorem 1.2].

Since $\widehat {P^*}(i)=\hat P_0(i,d-i)=\hat P(i,d-i)\not=0$ for some $i\in[0,k_1]$,
we have $m_{P^*}(0)\le k_1$. Similarly $\widehat {P^*}(d-j)\not=0$
for some $j\in[0,k_2]$ and hence $\deg P^*(x)\ge d-k_2$.

Set $f^*(x)=f_0(x,1)=P^*(x)(x+1)^{|C|}$.
Recall that $\widehat {f^*}(\kappa)=\hat f(\kappa,k_1+k_2-\da-\kappa)=0$
for all $\kappa\in[k_1-\da,k_1]$. Since $x^{k_1+1}\nmid f^*(x)$
and $\deg f^*(x)=|C|+\deg P^*(x)\ge |C|+d-k_2=k_1-\da$, by Lemma 1
there exists a $q\in\Cal P(p)$ such that
$N_q(f^*)\ge (k_1+1)-(k_1-\da-1)=\da+2$.

If $m_{f^*}(-1)=m_{P^*}(-1)+|C|<p$, then $N(P^*)=N(f^*)\ge
N_q(f^*)-1\ge k_1+k_2-d-|C|+1$,
therefore
$$|C|\ge k_1+k_2+1-d-N(P^*)=|A|+|B|-1-d-N(P^*).$$
This concludes our proof. \qed

\Ack. The authors are indebted to the two referees for their
many helpful comments. The revision was done during the second
author's visit to the Institute of Mathematics, Academia Sinica
(Taiwan); he would like to thank the Institute for its financial support.

\enddocument